\documentclass[a4paper,11pt]{article} % utilise cpde.cls
\usepackage{amsmath,amssymb,latexsym}
\catcode`\é = \active\def é{\'e} \catcode`\è = \active\def è{\`e}
\catcode`\ç = \active\def ç{\c{c}} \catcode`\à = \active\def
à{\`a} \catcode`\ê = \active \def ê{\^e} \catcode`\ù = \active
\def ù{\`u} \catcode`\ô = \active \def ô{\^o} \catcode`\û =
\active \def û{\^u} \catcode`\î = \active \def î{\^{\i}}
\catcode`\â = \active \def â{\^a} \catcode`\ö = \active \def
ö{{\"o}}

\def\Section#1{\section{#1}}

\newtheorem{thm}{Theorem}[section]  \newtheorem{cor}[thm]{Corollary}
\newtheorem{lem}[thm]{Lemma}   
  \newtheorem{prop}[thm]{Proposition}

\def\remark{\refstepcounter{thm}\bigskip\noindent\bf Remark \thethm\rm\ }
\newcommand{\preuve}[1][\!\!]{\bigskip\noindent{\bf Proof #1. \ \ }}
\def\fin{\hfill$\Box$\\}

\def\aaa{{\cal A}}\def\cc{{\cal C}}\def\dd{{\cal D}}
\def\hh{{\cal H}}\def\lll{{\cal L}}
\def\kk{{\cal K}}
\def\mm{{\cal M}}
\def\ss{{\cal S}}

\def\R{\mathbb R}\def\C{\mathbb C}
     
\def\D{\partial}\def\eps{\varepsilon}\def\phi{\varphi}
\def\norm#1{\left\Vert#1\right\Vert}

\def\set#1{\left\{#1\right\}}

\def\sep#1{\left(#1\right)}
\def\Re{{\mathrm Re\,}} 
 
\def\defegal{\stackrel{\text{\rm def}}{=}}

\def\Hess{\text{\rm{Hess}}}
\def\Id{\textrm{Id}}
\def\Ker{\textrm{Ker}}

\def\rhoinf{\rho_\infty}
\def\muinf{\mu_\infty}

\title{Hypocoercivity and exponential time decay for the linear
inhomogeneous relaxation Boltzmann equation}
\author{Fr{\'e}d{\'e}ric H{\'e}rau \footnote{\coord} \\
 Universit{\'e} de Reims }
\def\coord{Laboratoire de Mathématiques, UFR Sciences exactes et
naturelles, Université de Reims Champagne-Ardenne, Moulin de la
Housse, BP 1039 51687 Reims cedex 9, \tt herau@univ-reims.fr \rm }
\date{March 16, 2005}

\begin{document}
 \maketitle \thispagestyle{empty}
 \bibliographystyle{plain}
 \begin{abstract}
   We consider an inhomogeneous linear Boltzmann equation,
    with an external
    confining potential. The collision operator is
    a simple relaxation toward
     a local Maxwellian, therefore without diffusion.
     We prove the exponential time decay toward the global Maxwellian,
     with an explicit rate of decay. The methods are based
     on hypoelliptic methods transposed here to get spectral information.
     They were inspired by former works on the Fokker-Planck equation
     and the main feature of this work is that they are relevant
     although the equation itself has no regularizing properties.
 \end{abstract}

\section{Introduction.} This article is
devoted to the study of the long time behavior  of the solutions
of the following kinetic equation  in $\R_{t,x,v}^{1+2d}$ of
unknown $f$
\begin{equation}
  \label{LIB}
  \left\{
    \begin{array}[c]{l}
     \D_{t}f+v.\D_{x}f-\D_{x}V(x).\D_{v}f = Q(f),\\
    f|_{t=0}=f_0.
    \end{array}
\right.
\end{equation}
The right-hand side is a simple linear model for the Boltzmann
operator
\begin{equation}
Q(f) = \gamma( \rho \muinf - f),\ \ \ \  \rho(t,x) = \int
f(t,x,v)dv,
\end{equation}
where $\muinf$ is the Maxwellian in the velocity direction
\begin{equation}
\muinf(v) = \frac{e^{-v^2/2}}{(2\pi)^{d/2}}.
\end{equation}
This equation describes a system of  large number of particles
submitted to an external force deriving from a potential $V(x)$,
and for which the collision operator in the right-hand side  is a
simple relaxation toward the local Maxwellian $\rho\muinf$. In
particular there is no diffusion.

We suppose  that the derivatives of $V$  of order 2 or more are
bounded, and also that   $e^{-V} \in L^1$, which implies that
there is a unique steady state. In this case we say that $V$ is a
confining potential (anyway the adaptation in the case when
$e^{-V} \not\in L^1$ is straightforward, see remark \ref{notl1}).
It can be useful
 introduce the spatial Maxwellian and global
Maxwellian defined respectively by
$$
\rhoinf(x) = \frac{e^{-V(x)}}{\int e^{-V(x)} dx}, \ \ \  \mm(x,v)
= \rhoinf(x) \muinf(v).
$$
All steady states in $\mathcal{S}'(\R^{2d})$ are proportional to
the Maxwellian $\mm$.  In order to study the exponential decay we
now introduce an additionnal operator,
\begin{equation}
\Lambda^2 =-\gamma \D_v(\D_v +  v) -\gamma \D_x(\D_x + \D_x V) +
1.
\end{equation}
This operator has nice properties in the following weighted space
$$
B^2 = \set{ f \in \dd' \text{ s. t. } f/\mm^{1/2} \in L^2( dxdv)},
$$
with the natural norm defined by $\norm{f}_{B^2}^2=\int |f|^2
\mm^{-1}  dxdv$. Indeed  the closure from $\cc^\infty_0$ of
$\Lambda^2-1$   in $B^2$ is maximal accretive (see \cite{HelN04})
and has $0$ as single eigenvalue associated with the eigenfunction
$\mm$. We shall assume  the following:
\begin{equation} \label{spect}
 \textrm{Operator $\Lambda^2-1$ has a spectral gap $\alpha>0$ in
$B^2$}.
\end{equation}
Recall that the spectral gap is defined has the infimum of the
spectrum except the lowest eigenvalue.  We mention now some simple
cases when it happens. For example when $\Hess V \geq \lambda \Id$
then $\alpha=\lambda$. It is a special case of the one when
$|V'(x)| $ goes to infinity with $x$, which implies that
$\Lambda^2$ is with compact resolvent in $B^2$ and that
(\ref{spect}) is also satisfied. We refer to \cite{HN04} or
\cite{HelN04} and reference therein for complementary information
about it.

Now about the collision operator $Q$, we just mention here that it
is mass and positivity preserving and "dissipative" in the sense
non-negative in $B^2$ (see \cite{CCG03}).  We shall study it more
carefully later, and   refer to remark \ref{needQ} here for
complements.

It is easy to verify  that $-v\D_x+ \D_x V(x) \D_v +Q$ is also
mass and positivity preserving and dissipative, and that its
closure in $B^2$ from $\cc_0^\infty$ generates a semi-group of
contraction in $B^2$. The Cauchy problem (\ref{LIB}) is therefore
well posed and it was proven in \cite{CCG03}  that under
regularity assumptions and bounds on the solution $f$ of
(\ref{LIB}), $f(t)$ tends to $ \mm$ when $t$ goes to infinity
faster than any inverse power of $t$. We now state our main
result:

\begin{thm} \label{main}
There exists a constant $A >0$ depending only on the second and
third order derivatives of $V$, such that for all $L^1$ normalized
function $f_0 \in B^2$, we have the following
$$
\norm{f(t,.) - f_\infty}_{B^2} \leq 3 \norm{f_0- f_\infty}_{B^2}
e^{-\alpha^2 t /A}
$$
here $f_\infty = \mm$, and $f$ is the unique solution of equation
(\ref{LIB}).
\end{thm}

As a direct consequence we also obtain  the decrease of the
so-called relative entropy:

\begin{cor} \label{logdecay}
Under the hypothesis of the preceding Theorem, and assuming in
addition that  $f_0 \geq 0$, we have
$$
0\leq H(f,f_\infty)(t)  \defegal \iint f(t) \ln
\sep{\frac{f(t)}{f_\infty}}dxdv \leq 3 \norm{f_0}_{B^2} \norm{ f_0
- f_\infty}_{B^2}  e^{-\alpha^2 t/A}.
$$
\end{cor}

This study is motivated by proving the validity of some new tools
and ideas, namely the one called hypocoercivity, appeared in a few
recent articles in order to prove exponential time-decay
convergence for some inhomogeneous (mostly linear) kinetic
equations such as Fokker-Planck
\cite{HN04}\cite{HelN04}\cite{Vil05}\cite{HSS05},
Vlasov-Poisson-Fokker-Planck \cite{HVPFP05}, chains of anharmonic
oscillators \cite{Vil05}.  It is well-known that in the
homogeneous case the exponential decay can be easily obtained by
spectral methods if we assume (in the linear case) some coercivity
of the collision operator. In the inhomogeneous case, the global
coercivity is false in general, but can be obtained in a  modified
but norm-equivalent Hilbert space
 ($B^2$ in this work). This property can serve as a definition
 of hypocoercivity. Several tools can be used for this,
 essentially inspired by hypoelliptic ideas, that's why this name
 was introduced very recently in \cite{Vil05}. We mention some of
 them: The use of Kohn's method to get simultaneously
 hypocoercivity and hypoellipticity \cite{HN04}, \cite{HelN04},
\cite{Vil05}, in the case of the Fokker-Planck operator; the use
of analytic dilation and complex FBI-Bargmann transform
\cite{HSS05},  a method of multiplier via pseudodifferential
operators \cite{HSS05} or functional analysis using harmonic
oscillators and Witten Laplacian (linear part of \cite{HVPFP05}).

About the trend to equilibrium, this has been studied for the long
time, and we only want to quote the  entropy dissipating methods
introduced by Villani and Desvillettes to prove arbitrary and
explicit algebraic time decay. It was used for the Fokker-Planck
equation \cite{DV01}, for the model studied here \cite{CCG03}, and
in its main achievement for the full Boltzmann equation
\cite{DV04}. Now the question naturally arose whether the
exponential decay, obtained via hypoelliptic tools is also true
for non-hypoelliptic operators. In this work we choose a simple
example
 of collision operator which has no regularity property. It
 appears that Lie techniques (those also in the core of the
 hypoelliptic theory) also give sufficient information on the
 spectrum and in particular the spectral gap (in some modified
 $L^2$ space) to get hypocoercivity an then exponential decay.
  We therefore hope this techniques to be applied in
 the future to
 other inhomogeneous kinetic equations with linear or non-linear
  collision operators with or
 without regularity
 properties (see the review \cite{Vil03} for examples).

Eventually our result answer a question raised by C\'aceres,
Carillo and Goudon in  \cite{CCG03} about the applicability of
hypoelliptic techniques of  \cite{HN04} to obtain explicit
exponential decay of the model studied here.

\tableofcontents

%%%%%%%%%%%%%%%%%%%%%%%%%%%%%%%%%%%%%%%%%%%%%%%%%%%%%%%%%%%%%%%%%
%%%%%%%%%%%%%%%%%%%%%%%%%%%%%%%%%%%%%%%%%%%%%%%%%%%%%%%%%%%%%%%%%
%%%%%%%%%%%%%%%%%%%%%%%%%%%%%%%%%%%%%%%%%%%%%%%%%%%%%%%%%%%%%%%%%
%%%%%%%%%%%%%%%%%%%%%%%%%%%%%%%%%%%%%%%%%%%%%%%%%%%%%%%%%%%%%%%%%
%%%%%%%%%%%%%%%%%%%%%%%%%%%%%%%%%%%%%%%%%%%%%%%%%%%%%%%%%%%%%%%%%
%%%%%%%%%%%%%%%%%%%%%%%%%%%%%%%%%%%%%%%%%%%%%%%%%%%%%%%%%%%%%%%%%
%%%%%%%%%%%%%%%%%%%%%%%%%%%%%%%%%%%%%%%%%%%%%%%%%%%%%%%%%%%%%%%%%

\Section{Notations and functional analysis}

 For a Hilbert-space type analysis, it is more
convenient to work directly on the  equation  obtained after
conjugating with $\mm^{1/2}$. Therefore we pose $u=f/\mm^{1/2}$
which is now supposed to vary in $L^2$,  and the equation
satisfied by $u$ is then
\begin{equation} \label{eq.evolution}
 \left\{
    \begin{array}[c]{l}
     \D_{t}u+Ku = 0\\
    u|_{t=0}= u_0,
    \end{array}
\right.
\end{equation}
where we introduced the notations
\begin{equation} \label{defK}
 \left\{
    \begin{array}[c]{l}
     K= X_0 + \gamma(\Id -\Pi_1), \\
     X_0 = v.\D_{x}-\D_{x}V(x).\D_{v}, \\
     \Pi_1 u(x,v) = \sep{\int u(x,v') \muinf^{1/2}(v') dv'}
    \muinf^{1/2}(v) \ \ \ \textrm{for } u \in L^2.
    \end{array}
\right.
\end{equation}
It is immediate to see that $\Pi_1$ is an orthogonal projector
onto
$$
 E_1 = \set{u \textrm{ s.t. exists } \rho \in L^2(dx)
\textrm{ with } u= \rho \muinf^{1/2} }
$$
which is a closed subspace  of $L^2$. We therefore directly get
that in the new formulation in $L^2$, the collision operator,
which is now $-(\Id -\Pi_1)$ is dissipative.
 Of course we also have
$$
P_1 \mm^{1/2} = \mm^{1/2},
$$
 so that the square root of Maxwellian is in the kernel of the (new) collision
 operator. It is clear that $K$  and its adjoint
 $K^*=-X_0 + \gamma(\Id -\Pi_1)$ are
 well defined in $\mathcal{S}(\R^{2d})$,
in $\mathcal{S}'(\R^{2d})$ and as (non closed) operators in
$L^{2}(\R^{2d})$
 with domain $D(K)=D(K^{*})=\mathcal{S}(\R^{2d})$. We denote by the
 same later their maximal closure in $L^2$.
For $j=1,\ldots,d$, we introduce the differential operators
$a_{j}$, the annihilation operator $b_{j}$:
$$
a_{j}=\gamma^{1/2}\left(\D_{x_{j}}+\D_{x_{j}}V(x)/2\right)\qquad
b_{j}=\gamma^{1/2}\left(\D_{v_{j}}+v_{j}/2\right),
$$
and their formal adjoints
$$
a_{j}^{*}=\gamma^{1/2}\left(-\D_{x_{j}}+\D_{x_{j}}V(x)/2\right)\qquad
b_{j}^{*}=\gamma^{1/2}\left(-\D_{v_{j}}+v_{j}/2\right).
$$
For the sake of conciseness, the letters $a$ and $b$ denote the
vectors
$$
a=\left(
  \begin{array}[c]{l}
   a_{1}\\
   \vdots\\
   a_{d}
  \end{array}
\right) \qquad b=\left(
  \begin{array}[c]{l}
   b_{1}\\
   \vdots\\
   b_{d}
  \end{array}
\right)
$$
while $a^{*}$ and $b^{*}$ are the forms
$$
a^{*}=\left(a_{1}^{*},\ldots,a_{d}^{*}\right) \quad
b^{*}=\left(b_{1}^{*},\ldots,b_{d}^{*}\right).
$$
Up to the factor $\gamma$, the non negative operator
$\Lambda^{2}-1$ is nothing but the sum of the Witten Laplacian on
$0$-forms $
\gamma^{-1}a^{*}a=-\Delta_{x}+\left|\D_{x}V(x)\right|^{2}/4-
\Delta V(x)/2 $ and of the harmonic oscillator $
\gamma^{-1}b^{*}b=-\Delta_{v}+v^{2}/{4}-d/2. $ Under the
hypothesis on $V$, $\Lambda^2-1$  is a nonnegative,
$\mathcal{S}(\R^{2d})$ is a core and $\Lambda^{r}$ is well-defined
for $r\in\R$. Moreover the kernel of $\Lambda^{2}-1$ is $E_0 = \C
\mm^{1/2}$. For the following we shall denote $\Pi_0$ the
orthogonal projector onto $E_0$,
$$
\Pi_{0}u= (u,\mm^{1/2})_{L^2} \mm^{1/2}
$$
which also extends to $u\in \mathcal{S}'(\R^{d})$. Let us now
recall some relations involving the operators $a$, $b$ and $X_0$,
that can be found for example in section 1 of  \cite{HN04}. We
mention here that since the operators are continuous in $\ss$ and
$\ss'$ there is no problem of defining their commutators. First
recall the Canonical Commutation Relations
 (CCR) for $b$ and their counterparts for $a$
  \begin{equation}
    [b_{j},b_{k}]=[b_{j}^{*},b_{k}^{*}]=0
    \qquad
    [b_{j},b_{k}^{*}]=\gamma\delta_{jk}
    \qquad
    [a_{j},a_{k}]=[a_{k},a_{j}]=0
    \qquad
    [a_{j},a_{k}^{*}]=\gamma\D^{2}_{x_{j}x_{k}}V.
  \end{equation}
It is also clear that the  $a$'s and $a^*$'s commute with the
$b$'s and $b^*$'s. The main remark is that the $a_{j}$'s,
$a_{j}^{*}$'s are in the Lie algebra generated by the $b_{j}$'s,
$b_{j}^{*}$'s and the vector field $X_{0}$:
\begin{equation}
  \label{eq.LieX0b}
  [b_{j},X_{0}]=a_{j}\qquad [b_{j}^{*},X_{0}]=a_{j}^{*}.
\end{equation}
 Similarly, the $b_{j}$'s and $b_{j}^{*}$'s can be derived from
the $a_{j}$'s, $a_{j}^{*}$'s and $X_{0}$
\begin{equation}
  \label{eq.LiemodX0a}
  [a_{j},X_{0}]=-\sum_{k=1}^{d}\left(\D^{2}_{x_{j}x_{k}}V\right)b_{k}
\qquad
  [a_{j}^{*},X_{0}]
=-\sum_{k=1}^{d}b_{k}^{*}\left(\D^{2}_{x_{k}x_{j}}V\right).
\end{equation}
The relations (\ref{eq.LieX0b}) and (\ref{eq.LiemodX0a}) are
summarized by
\begin{equation}
  \label{eq.abcompacte}
  [b,X_{0}]=a\ ,\
[b^{*},X_{0}]=a^{*}\ ,\ [a,X_{0}]=-\Hess V b \ \text{and}\
[a^{*},X_{0}]=- b^{*}\Hess V.
\end{equation}
By combination we have the useful formulas:
 \begin{eqnarray}
    \label{commm}
    &&
    \left[\Lambda^{2},X_{0}\right]= -b^{*}(\Hess V -\Id)a
    -a^{*}(\Hess V -\Id)b,\\
    &&
    b^{*}(a^{*}a)=(a^{*}a)b^{*}\qquad a^{*}(a^{*}a)=(a^{*}a)a^{*}-\gamma
    a^{*}\Hess V_{\beta}\\
    \text{and}&&
    a^{*}(b^{*}b)=(b^{*}b)a^{*}\qquad b^{*}(b^{*}b)=(b^{*}b)b^{*}-\gamma
    b^{*},
  \end{eqnarray}
and their adjoint relations hold as equality of continuous
operators in $\mathcal{S}(\R^{2d})$ and $\mathcal{S}'(\R^{2d})$.

%%%%%%%%%%%%%%%%%%%%%%%%%%%%%%%%%%%%%%%%%%%%%%%%%%%%%%%%%
%%%%%%%%%%%%%%%%%%%%%%%%%%%%%%%%%%%%%%%%%%%%%%%%%%%%%%%%%%%%
%%%%%%%%%%%%%%%%%%%%%%%%%%%%%%%%%%%%%%%%%%%%%%%%%%%%%%%%%%%%
%%%%%%%%%%%%%%%%%%%%%%%%%%%%%%%%%%%%%%%%%%%%%%%%%%%%%%%%%

\Section{Hypocoercivity for operator $K$.} \label{se.isohyp}

In this section we continue to work with operator $K$ defined in
the preceding section. We shall prove that operator $K$ is
hypocoercitive, i.e. coercitive in $L^2$ with a modified scalar
product. For this we introduce an additional operator
$$
L= \Lambda^{-2}a^* b = \Lambda^{-2}(\sum_j a^*_j b_j)
$$
We shall see later that this operator is explicitly bounded in
terms of the second and third derivatives of $V$.

\begin{prop} \label{hypoco} Let $\alpha$ be defined in
(\ref{spect}).
  Then there exists $\eps, A >0$ such that for all
   $\ss \ni u \perp \mm^{1/2}$ we have
$$
\Re \sep{ Ku, (\Id + \eps(L+L^*))u} \geq  \frac{\alpha^2}{A}
\norm{u}^2,
$$
where  $A$ can be chosen to depend explicitly on $\gamma$ and  the
second and third derivatives of $V$, and $\norm{\eps L} \leq 1$.
\end{prop}

\preuve Let us take $u\in L^2$ and $\eps>0$. We write
\begin{equation*}
\begin{split}
& \Re \sep{ Ku, (\Id + \eps(L+L^*))u} \\
&= \Re \gamma \sep{ (\Id -\Pi_1)u, (\Id + \eps(L+L^*))u} +\Re
\sep{ X_0 u,
(\Id + \eps(L+L^*))u} \\
&= \gamma \norm{ (Id-\Pi_1) u}^2 + \eps\gamma \Re \sep{ (\Id
-\Pi_1)u,
(L+L^*)u} +\eps \Re \sep{ X_0 u,+ (L+L^*)u}\\
&= I + II + III,
\end{split}
\end{equation*}
where in the last term we used the fact that $X_0$ is
skew-adjoint. We first study the two first terms. Using the
Cauchy-Schwartz inequality we can write
\begin{equation} \label{partauto}
I + II \geq \frac{\gamma}{2}\norm{ (Id-\Pi_1) u}^2 -
\frac{\gamma}{2}\eps^2 \norm{ (L+L^*)u}^2 \geq
\frac{\gamma}{2}\norm{ (Id-\Pi_1) u}^2 - \eps^2 \gamma \norm{L}^2
\norm{u}^2.
\end{equation}
Now we study more carefully the third one
\begin{equation} \label{partauto1}
III= \eps \Re \sep{ X_0 u, (L+L^*)u} = \eps \Re \sep{ [L,X_0]
u,u},
\end{equation}
again since $X_0$ is skewadjoint.  Recalling that $L=
\Lambda^{-2}a^* b$  we can write using the equalities
(\ref{commm})
\begin{equation} \label{partauto2}
\begin{split}
 [L,X_0] & = [\Lambda^{-2} a^*b, X_0] \\
 & =[\Lambda^{-2}, X_0] a^* b + \Lambda^{-2}[a^*,X_0]b
 + \Lambda^{-2}a^*[b,X_0]\\
 & =-\Lambda^{-2}[\Lambda^2,X_0] \Lambda^{-2}a^*b - \Lambda^{-2}b^*\Hess V b
  + \Lambda^{-2}a^*a.
 \end{split}
 \end{equation}
  Here we used the fact that for $A$, $B$, and $B^{-1}$ continuous on
  $\ss$ and $\ss'$ we have $[A, B^{-1}] = -B^{-1}[A,B]B^{-1}$.
Let us denote
$$
\aaa = -\Lambda^{-2}[\Lambda^2,X_0] \Lambda^{-2}a^*b -
\Lambda^{-2}b^*\Hess V b
$$
We postpone to the end of this section the proof of following
lemma concerning operators $\aaa$ and $L$.
\begin{lem} \label{aaa}  Operator $\aaa$ and $L$ are
 bounded on $L^2$. Besides their norms of  can be explicitly bounded in
terms of $\gamma$ and the second and third  order derivatives of
$V$.
\end{lem}
Now it is clear that denoting $\hh^{0,-1} = \set{ bu \textrm{ for
} u\in L^2}$, operator $b$ considered as an operator from $L^2$
into $\hh^{0,-1}$ satisfies $E_1 \subset \Ker (b)$. In fact this
comes from the fact that $b$ is then the annihilation operator in
the velocity variable. Since it appears only in the right in the
expression of $\aaa$ we therefore get
$$
E_1 \subset \Ker( \aaa).
$$
As a consequence we can write $\aaa = \aaa (I-\Pi_1)$ and
therefore
$$
[L,X_0] = \aaa(\Id - \Pi_1) + \Lambda^{-2}a^*a
$$
Putting this in the expression of the term III appearing in
(\ref{partauto1}) yields
\begin{equation}
\begin{split}
III & = \eps \Re \sep{ \aaa (\Id - \Pi_1) u,u} + \eps \Re \sep{
\Lambda^{-2}a^*a u,u} \\
    & \geq -\frac{\gamma}{4}\norm{(\Id - \Pi_1) u}^2
    -\frac{\eps^2}{\gamma}\norm{\aaa}^2\norm{u}^2 +\eps \Re \sep{
\Lambda^{-2}a^*a u,u}.
\end{split}
\end{equation}
Now it is also clear that $\Lambda^2$,  $a^*a$ and $P_1$  commute
together, and we can therefore write for the second term appearing
in the last inequality
\begin{equation} \label{bornepi}
\begin{split}
\eps \Re \sep{ \Lambda^{-2}a^*a u,u} &=\eps \Re \sep{
\Lambda^{-2}a^*a \Pi_1 u,u} + \eps \Re \sep{ \Lambda^{-2}a^*a
(1-\Pi_1) u,u} \\
& =\eps \Re \sep{ \Lambda^{-2}a^*a \Pi_1 u, \Pi_1 u} + \eps \Re
\sep{ \Lambda^{-2}a^*a
(1-\Pi_1) u,(1-\Pi_1)u} \\
&  \geq \eps \Re \sep{ \Lambda^{-2}a^*a \Pi_1 u, \Pi_1 u} - \eps
\norm{(1-\Pi_1)u}^2.
\end{split}
\end{equation}
for the last inequality we simply used the fact that $a^*a \leq
\Lambda^2$, which implies easily that the norm of $a\Lambda^{-1} $
is bounded by $1$.

Now we can take into account the spectral gap property of
$\Lambda^2$ together with the fact that $u$ is supposed to
orthogonal to $\mm^{1/2}$. We write it as a lemma

\begin{lem} Recalling that $u \perp \mm^{1/2}$ we have
$
\Re \sep{ \Lambda^{-2}a^*a \Pi_1 u, \Pi_1 u} \geq \frac{\alpha}{1
+\gamma} \norm{w}^2
$
\end{lem}
\preuve We first notice that $\Lambda^2$ and $a^*a$ leave $E_1$
invariant, and that
$$
\Lambda^2|_{E_1} = 1 + a^*a
$$
is essentially the Witten Laplacian in the spacial direction. We
define now $\tau$ to be the spectral gap for $a^*a$. Now
 for  $w \in \ss(\R^d_x)$ such that $w \perp \rhoinf^{1/2}$ we
have
$$
(a^*a w,w)_{L^2(\R^d_x)} \geq \tau \norm{w}^2_{L^2(\R^d_x)}.
$$
Since $(a^*a + 1)^{-1/2}w \perp \rhoinf^{1/2}$ also  we get
$$
(a^*a (a^*a + 1)^{-1/2}w,(a^*a + 1)^{-1/2}w)_{L^2(\R^d_x)} \geq
\frac{\tau}{1 + \tau} \norm{w}^2_{L^2(\R^d_x)}.
$$
Now the following inequalities are clear:
$$
\frac{\tau}{1 + \tau} \geq \frac{\alpha}{1 + \alpha} \geq
\frac{\alpha}{1 + \gamma}.
$$
Indeed $\tau \geq \alpha$ from the definition of $\alpha$ and
$\alpha \leq \gamma$ because of the harmonic part of
$\Lambda^2-1$. Now since $\Pi_1 u \in E_1 \cap E_0^\perp$ we get
the result of the lemma from the preceding study by applying it to
the function defined for a.e. $v$ by $x \mapsto \Pi_1 u(x,v)$.
\fin

\noindent \bf End of the proof of Proposition \ref{hypoco}.\rm \ \
Now we can put  the result of the lemma in (\ref{bornepi}) and we
get
\begin{equation} \label{bornepi2}
\begin{split}
\eps \Re \sep{ \Lambda^{-2}a^*a u,u} &  \geq \eps
\frac{\alpha}{1+\gamma} \norm{\Pi_1 u}^2 - \eps
\norm{(1-\Pi_1)u}^2.
\end{split}
\end{equation}
We obtain the following lower bound for the term III from
(\ref{partauto})
\begin{equation}
\begin{split}
III &  \geq -\frac{\gamma}{4}\norm{(\Id - \Pi_1) u}^2
    -\frac{\eps^2}{\gamma}\norm{\aaa}^2\norm{u}^2 + \eps
\frac{\alpha}{1+\gamma} \norm{\Pi_1 u}^2 - \eps
\norm{(1-\Pi_1)u}^2.
\end{split}
\end{equation}
Eventually putting together the estimate on I+II and III we get
\begin{equation*}
\begin{split}
& \Re \sep{ Ku, (\Id + \eps(L+L^*))u} \\
& \geq  \frac{\gamma}{2}\norm{ (Id-\Pi_1) u}^2 - \eps^2 \gamma
\norm{L}^2 \norm{u}^2
\\&  -\frac{\gamma}{4}\norm{(\Id - \Pi_1) u}^2
    -\frac{\eps^2}{\gamma}\norm{\aaa}^2\norm{u}^2 + \eps
\frac{\alpha}{\delta^2 + \alpha} \norm{\Pi_1 u}^2 - \eps
\norm{(1-\Pi_1)u}^2. \\
& \geq \frac{\gamma}{8}\norm{(\Id - \Pi_1) u}^2 + \eps
\frac{\alpha}{1+\gamma} \norm{\Pi_1 u}^2
-\eps^2\sep{\gamma^{-1}\norm{\aaa}^2+ \gamma \norm{L}^2}
\norm{u}^2
\end{split}
\end{equation*}
by taking $\eps \leq \gamma /8$.  Now we use the fact that $\Pi_1$
is an orthogonal projector and that $\eps \leq \gamma/8$:
$$
\Re \sep{ Ku, (\Id + \eps(L+L^*))u}  \geq \sep{ \eps
\frac{\alpha}{1+\gamma} - \eps^2\sep{\gamma^{-1}\norm{\aaa}^2+
\gamma \norm{L}^2}} \norm{u}^2.
$$
For $\eps/\alpha$ sufficiently small, but depending only on
$\gamma$ and the second and third order derivatives of $V$ via
lemma \ref{aaa}, we obtain
\begin{equation*}
\begin{split}
 \Re \sep{ Ku, (\Id + \eps(L+L^*))u} \geq \frac{\alpha^2}{A} \norm{u}^2
\end{split}
\end{equation*}
where $A$ satisfies the hypothesis of the Proposition.  The proof
is then complete. \fin

\preuve[ of lemma \ref{aaa}] Recall that $L= \Lambda^{-2} a^*b$.
Now for $\aaa$ we use the expression of $[\Lambda^2, X_0]$ in
(\ref{commm}) and we get
\begin{equation}
\begin{split}
\aaa = &  \Lambda^{-2} b^* (\Hess(V) - \Id) a \Lambda^{-2} a^* b +
    \Lambda^{-2} a^* (\Hess(V) - \Id) b \Lambda^{-2} a^* b
     - \Lambda^{-2} b^* \Hess(V)b.
\end{split}
\end{equation}
We therefore see that it is sufficient to prove that for any
$d\times d$ real matrix $M(x)$ depending only on $x$,   bounded
and with first derivative bounded, the following operators
$$
\Lambda^{-2} b^* M(x) a, \ \ \ \Lambda^{-2} b^* M(x) b, \ \ \
\Lambda^{-2} a^* M(x) b
$$
are bounded as operators on $L^2$.  We give the proof for the
first one since for the remaining ones, the proof is similar and
easier.
 We shall prove the result for its adjoint
$a^* M(x) b \Lambda^{-2}$. We write for  $u\in \ss$,
\begin{equation} \label{detail}
\begin{split}
\norm{a^* M(x) b \Lambda^{-2}u}  & \leq \sum_{j,k} \norm{a_j^*
M_{j,k}(x) b_k \Lambda^{-2}u} \\
& \leq  \sum_{j,k} \norm{M_{j,k}(x) a_j^* b_k \Lambda^{-2}u} +
\sum_{j,k} \gamma^{1/2} \norm{(\D_{x_j} M_{j,k})(x) b_k
\Lambda^{-2}u} \\
& \leq \sep{ \norm{M}_{L^{\infty}} + \gamma^{1/2} \norm{\D_x
M}_{L^{\infty}}} \sum_{j,k} \sep{\norm{a_j^* b_k \Lambda^{-2}u} +
\norm{ b_k \Lambda^{-2}u}}
\end{split}
\end{equation}
where we used the fact that $[a_j^*, M_{j,k}] =-\D_{x_j} M_{j,k}$.
Now this is straightforward to check that $\norm{ b_k
\Lambda^{-2}u} \leq \norm{u}$ since $b^*b \leq \Lambda^2$ and $1
\leq \Lambda^2$. For the other term $\norm{a_j^* b_k
\Lambda^{-2}u}$ in the last inequality of (\ref{detail}), we write
\begin{equation} \label{detail2}
\begin{split}
& \norm{a_j^* b_k \Lambda^{-2}u}^2 = \sep{ a_ja_j^*b_k
\Lambda^{-2}u, b_k \Lambda^{-2}u} \\
& \leq \sep{ a_j^*a_j b_k \Lambda^{-2}u, b_k \Lambda^{-2}u} +
\sep{ \gamma(\D^2_{x_j} V) b_k \Lambda^{-2}u, b_k \Lambda^{-2}u} \\
& \leq \sep{ \Lambda^2 b_k \Lambda^{-2}u, b_k \Lambda^{-2}u} +
 \gamma\norm{ \Hess V}_{L^{\infty}}
\sep{  b_k \Lambda^{-2}u, b_k \Lambda^{-2}u} \\
& \leq \sep{ \Lambda^2 b_k \Lambda^{-2}u, b_k \Lambda^{-2}u} +
 \gamma \norm{ \Hess V}_{L^{\infty}} \sep{  b_k \Lambda^{-2}u, b_k
\Lambda^{-2}u}
\end{split}
\end{equation}
Now using the fact that $[\Lambda^2, b_k] = -\gamma b_k$ we can
continue the preceding series of inequalities:
\begin{equation*}
\begin{split}
& \leq \sep{ \Lambda^2 b_k \Lambda^{-2}u, b_k \Lambda^{-2}u} +
 \gamma\norm{ \Hess V}_{L^{\infty}}
\sep{  b_k \Lambda^{-2}u, b_k \Lambda^{-2}u} \\
& \leq \sep{  b_k u, b_k \Lambda^{-2}u} + \gamma (\norm{ \Hess
V}_{L^{\infty}}+ 1)
\sep{  b_k \Lambda^{-2}u, b_k \Lambda^{-2}u} \\
& \leq  \gamma (\norm{ \Hess V}_{L^{\infty}}+ 2)
\norm{u}^2 \\
\end{split}
\end{equation*}
again since $b^*_kb_k \leq \Lambda^2$ and $1 \leq \Lambda^2$.
Therefore the  term $\norm{a_j^* b_k \Lambda^{-2}u}$ in the last
inequality of (\ref{detail}) is also bounded by $C\norm{u}$, where
$C$ depends only on $\gamma$  and the second  and third
derivatives of $V$. The proof of lemma \ref{aaa}
 is therefore complete. \fin

%%%%%%%%%%%%%%%%%%%%%%%%%%%%%%%%%%%%%%%%%%%%%%%%%%
%%%%%%%%%%%%%%%%%%%%%%%%%%%%%%%%%%%%%%%%%%%%%%%%%%%
\section{Proof of the Theorem and comments}

We go on studying operator $K$ defined in the preceding sections.
We first quote an easy result from \cite{HVPFP05} relying
hypocoercivity to exponential decay

 \begin{lem}[lemma A6 in \cite{HVPFP05}] \label{expdecayabs}
Let $\kk$ be the infinitesimal generator of a semigroup of
contraction on a Hilbert space $H$ and suppose that there exist a
constant  $\delta >0$ and a bounded operator  $\lll$ with norm
bounded by $C_\lll \geq 1$
 such that
\begin{equation} \label{abstr}
\forall u  \in D(\kk), \ \ \ \delta\norm{u}^2 \leq  \Re(\kk u,u) +
\Re(\kk u, (\lll+\lll^*)u)
\end{equation}
then for all $u_0 \in H$ and $t \geq 0$ we have $ \norm{ e^{-t\kk}
u_0} \leq 3 e^{-\frac{\delta t }{3C_\lll}} \norm{u_0}$
\end{lem}

\preuve[ of Theorem \ref{main}] we first consider operator $K$
defined in the preceding sections.The result for $K$ in the
Hilbert space ${(\mm^{1/2})}^{\perp}$ is then a direct consequence
of Proposition \ref{hypoco}. Indeed it suffices to apply the
preceding lemma  when
 considering  $K$ in place of $\kk$, $\eps L$
 in place of $\lll$, $C_\lll=1$ and $ \alpha^2/A$ in place of $\delta$.
 The result of the theorem is then a simple transcription in terms of
  $f= \mm^{1/2} u$ and $f-f_\infty = \mm^{1/2} (u-\Pi_0 u)$ (of course
  $(u-\Pi_0 u) \perp \mm^{1/2}$, which is an other way to say that
  $\iint (f-f_\infty)dxdv = 0$). The proof is complete. \fin

 \preuve[ of Corollary \ref{logdecay}]
This is then a direct consequence of the main Theorem, and the
proof follows exactly the one in \cite[corollary 0.2]{HN04}:  Let
$f_{0}$ be a $L^1$-normalized function which belongs to
$\mm^{1/2}L^2$ and let $f(t)$ be the solution of (\ref{LIB}), (it
 stays non-negative for all time). The non-negativity of the relative entropy
   is clear from the fact that
 $f$ stays $L^1$ normalized. For the other side,  we write for $t\geq 0$,
with $f_\infty= \mm$:
\begin{equation*}
 H(f(t)|\mm)=\iint f(t)\ln\left(\frac{f(t)}{\mm}\right)~dxdv
=\iint
\frac{f(t)}{\mm^{1/2}}\mm^{1/2}\log\left(\frac{f(t)}{\mm}\right)~dxdv.
\end{equation*}
Applying first  $\ln(x)\leq x-1$ and then the Cauchy-Schwarz
inequality for $t\geq 0$, we get
\begin{equation}
\begin{split}
   H(f(t)|\mm) & \leq \iint
   \frac{f(t)}{\mm^{1/2}}\mm^{1/2}\left(\frac{f(t)}{\mm}-1\right)~dxdv
   \\
   & \leq
   \left\|\frac{f(t)}{\mm^{1/2}}\right\|
   \left\|\frac{f(t)}{\mm^{1/2}}-\mm^{1/2}\right\|
    = \left\|f(t)\right\|_{B^2}
   \left\|f(t)- \mm\right\|_{B^2}.
 \end{split}
 \end{equation}
Hence the exponential decay of the relative entropy is a
consequence of Theorem \ref{main} . \fin

\bigskip We end this work by making some remarks about the main
result:

 \remark \label{notl1} In the case $e^{-V} \not\in L^1$,
Theorem \ref{main}  remains true when replacing  $f_\infty$ by
$0$. A careful study of the proof shows that it is exactly the
same with the following adaptations: $\rhoinf(x) = e^{-V(x)}$ is
not anymore in $L^1(dx)$ and neither does $\mm(x,v) = \rhoinf(x)
\muinf(v)$ in $L^1(dxdv)$. Then $\alpha$ is  the bottom of the
spectrum of operator $\Lambda^2-1$, and in the proof $E_0 =
\set{0}$ and the corresponding projector is  $\Pi_0 = 0$. In this
case Theorem \ref{main} has to be understood has a vanishing.

\remark \label{needQ} A natural question is to understand the
common features between for example the Fokker-Planck operator
(and also the chains of anaharmonic oscillators) and the one
studied here. Let us see this after conjugating by the square root
of the Maxwellian and only in $L^2$ for  the velocity variable
(homogeneous case). The collision
 operators are respectively
 $$
 Q_{fp} = b^*b \ \ \textrm{ (Fokker-Planck) },  \ \ \
  \ Q_{lib}= \Id- \Pi_1 \ \ \textrm{ (Linear Inhom. Boltzmann)}.
 $$
 Here $\Pi_1$ is the orthogonal projection on the space spanned by
 the  Maxwellian  in the
 velocity variable $\muinf^{1/2}$.
A simple remark can be made using the Hermite decomposition of
functions (in the velocity variable), which we denote  $H_k$.  We
know that the annihilation operators $b$ and its conjugate $b^*$
(creation operator) have both a nice description, since for all $k
\geq 1$,
$$
b H_k = \sqrt{k} H_{k-1}, \ \ \ b^* H_{k-1} = \sqrt{k} H_{k}
$$
and $b H_0 =0$ (recall that $H_0 = \muinf^{1/2}$). Therefore $b$
and $b^*$ can be represented by respectively an upperdiagonal and
a subdiagonal infinite matrices with coefficients $\sqrt{k}$.
Immediately we get the expression of the harmonic oscillator
$b^*b$, for which the Hermite polynomials are an Hilbertian base.
Now using this decomposition we can build operators $c$ and its
adjoint $c^*$ by imposing
$$
c H_k =  H_{k-1}, \ \ \ c^* H_{k-1} =  H_{k}
$$
and $c H_0 =0$. Then $c$ and $c^*$ have the same representation as
matrices than $b$ and $b^*$ respectively (note that they are
bounded). It is immediate to check that
$$
c^*c = \Id-\Pi_1 = Q_{lib}.
$$
As a conclusion, and transferring this in the inhomogeneous case,
we see that the Fokker-Planck and the linear inhomogeneous
Boltzmann models have the same structure, explaining (a bit) the
efficiency of hypoelliptic methods in the last case.

%%%%%%%%%%%%%%%%%%%%%%%%%%%%%%%%%%%%%%%%%%%%%%%%%%%%%%%%%%{\`u}
%%%%%%%%%%%%%%%%%%%%%%%%%%%%%%%%%%%%%%%%%%%%%%%%%%%%%%%%%%{\`u}
%%%%%%%%%%%%%%%%%%%%%%%%%%%%%%%%%%%%%%%%%%%%%%%%%%%%%%%%%%{\`u}
%%%%%%%%%%%%%%%%%%%%%%%%%%%%%%%%%%%%%%%%%%%%%%%%%%%%%%%%%%{\`u}
%%%%%%%%%%%%%%%%%%%%%%%%%%%%%%%%%%%%%%%%%%%%%%%%%%%%%%%%%%{\`u}
%%%%%%%%%%%%%%%%%%%%%%%%%%%%%%%%%%%%%%%%%%%%%%%%%%%%%%%%%%{\`u}
%%%%%%%%%%%%%%%%%%%%%%%%%%%%%%%%%%%%%%%%%%%%%%%%%%%%%%%%%%{\`u}

\end{document}